\documentclass[a4paper,12pt]{article}
\usepackage[english,russian]{babel}
\usepackage[T2A]{fontenc}
\usepackage[cp1251]{inputenc}
\usepackage[english]{babel}
\usepackage{amsthm}
\usepackage[tbtags]{amsmath}
\usepackage{amsfonts,amssymb}
\begin{document}

\newtheorem{theorem}{Theorem}
\newtheorem{lemma}{Lemma}
\newtheorem{proposition}{Proposition}
\newtheorem{Cor}{Corollary}

\begin{center}
{\large\bf Ideals and Factor Rings of Centrally Essential Rings}
\end{center}
\begin{center}
O.V. Lyubimtsev\footnote{Nizhny Novgorod State University; email: oleg lyubimcev@mail.ru .}, A.A. Tuganbaev\footnote{National Research University `MPEI', Lomonosov Moscow State University; email: tuganbaev@gmail.com .}
\end{center}

\textbf{Abstract.} It is proved that the ring $R$ with center $Z(R)$, such that the module $R_{Z(R)}$ is an essential extension of the module $Z(R)_{Z(R)}$, is not necessarily right quasi-invariant, i.e., maximal right ideals of the ring $R$ are not necessarily ideals. We use central essentiality to obtain conditions which are sufficient to the property that all maximal right ideals are ideals


The work of O. Lyubimtsev is done under the state assignment No~0729-2020-0055. A. Tuganbaev is supported by Russian Scientific Foundation, project 16-11-10013P.

\textbf{MSC2010 database 16D25, 16R99}

\section{Introduction}

We consider associative unital non-zero rings only. A ring $A$ is said to be \textsf{centrally essential} if either $A$ is commutative or for every non-central element $a\in A$, there exist such non-zero central elements $x,y$ that $ax = y$. It is clear that the ring $A$ with center $Z(A)$ is centrally essential if and only if the module $A_{Z(A)}$ is an essential extension of the module $Z(A)_{Z(A)}$. In addition, any commutative ring is centrally essential. 

For a ring $A$, we denote by $J(A)$ the Jacobson radical of $A$.

\textbf{1.1. Remark.} If $A$ is a ring and the factor ring $A/J(A)$ is centrally essential, then the ring $A/J(A)$ is commutative and ring $A$ is right and left quasi-invariant.

$\lhd$ Since $A/J(A)$ is a centrally essential semiprime ring, the ring $A/J(A)$ is commutative, by \cite[Proposition 3.3.]{MT18}. In particular, all maximal right (left) ideals of the ring $A/J(A)$ are ideals. Therefore, the ring $A$ is right and left quasi-invariant.~$\rhd$

\textbf{1.2. Remark.} All right semi-Artinian or semiperfect centrally essential rings are right quasi-invariant.

Remark 1.2 follows from Remark 1.1 and the property that the factor ring $A/J(A)$ centrally essential of the ring $A$ is commutative if the ring $A$ is semiperfect or right semi-Artinian; see \cite[Proposition 3.4]{MT19a} and \cite[Theorem 1.5]{MT19b}.

\textbf{1.3. Remark.} If $A$ is a centrally essential ring and the factor ring $A/J(A)$ is commutative, then any minimal right ideal $S$ of the ring $A$ is contained in $Z(A)$. In particular, all minimal right ideals of the ring $A$ are ideals and $\text{Soc}\,A_A$ is contained in the center of the ring $A$.

$\lhd$ The simple $A$-module $S_A$ is a cyclic module with non-zero generator $s\in S$. We set $K=S\cap Z(A)$. By assumption, there exist two non-zero central elements $x,y\in Z(A)$ such that $0\ne sx=y\in S$. Then $S=yA$ and $y\in K\neq 0$. Since $J(A)$ is the intersection of annihilators of all simple right $A$-modules, we have $SJ(A) = 0$. Since the ring $A/J(A)$ is commutative, we have $ab-ba\in J(A)$ for any two elements $a,b\in A$. In addition, $k(ab-ba)=0$ for any $k\in K$. Since $k\in Z(A)$, we have $(ka)b=kba= b(ka)$. This implies that $ka\in Z(A)$. Next, it follows from $ka\in S$ that $K$ is a right ideal of the ring $A$. Since $S$ is a minimal right ideal, we have $S=K\subseteq Z(A)$.~$\rhd$

In connection to Remarks 1.1, 1.2 and 1.3, we prove Theorem 1.4 which is the main result of this paper.

\textbf{1.4. Theorem.}\\
\textbf{a.} There exists a centrally essential ring containing a maximal right ideal which is not an ideal.

\textbf{b.} There exist a centrally essential ring containing a closed\footnote{A submodule $X$ of the module $M$ is said to be \textsf{closed} if $X=Y$ for any submodule $Y$ of $M$ which is an essential extension of the module $X$.} right ideal which is not an ideal. 

In connection to Theorem 1.4, we formulate an open question. 

\textbf{1.5. Open question.} Is it true that there exists a centrally essential ring with a minimal right ideal which is not an ideal? 

\section{The proof of Theorem 1.4}

\textbf{2.1. Lemma.} Let $A$ be a ring, $R=A[x]$ be the polynomial ring, and let $M$ be a maximal right ideal of the ring $A$.

\textbf{a.} $MR+xR$ is a maximal right ideal of the ring $R$.

\textbf{b.} If $MR+xR$ is an ideal of the ring $R$, then $M$ is an ideal of the ring $A$, the factor rings $A/M$ and $R/(MR+xR)$ are isomorphic division rings, $(A/M)[x]$ is a principal one-sided ideal domain, and we have an isomorphism $\alpha\colon R/(MR+xR)\to A/M$ such that $\alpha(f+(MR+xR))=f_0+M$, where $f_0$ is the constant term of the polynomial $f\in R$.

\textbf{c.} If $A$ is a division ring and the polynomial ring $R$ is right quasi-invariant, then the ring $A$ is a field.

\textbf{d; see \cite[Theorem 7]{HuhJKL02}.} If the polynomial ring $R$ is right quasi-invariant, then the factor ring $A/J(A)$ is commutative.

$\lhd$ \textbf{a.} Let $f=f_0+xg\in R$ be a polynomial which is not contained in the right ideal $MR+xR$ of the ring $R$, where $f_0\in A$ and $g\in R$. Then $f_0\notin M$, since otherwise $f=f_0+xg\in MR+xR$. Therefore, there exist elements $d\in A$ and $m'\in M$ such that 
$$
1=f_0d+m'=(f-xg)d+m'=fd-xgd+m'\in fR+xR+MR= fR+MR+xR.
$$
Therefore, $MR+xR$ is a maximal right ideal of $R$.

\textbf{b.} These well known assertions are directly verified.

\textbf{c.} Let $a$ and $b$ be two non-zero elements of the division ring $A$. Since $A$ is a division ring, the domain $R$ has the Euclidean algorithm. Therefore, $(a + x)R$ is a maximal right ideal of the right quasi-invariant domain $R$. Then $(a + x)R$ is an ideal of $R$. Therefore, $b(a+x)\in (a+x)R$ and $ba + bx = ac + xc$ for some $c\in A$, whence we have $bx=xc$ and $ba=ac$. Then $b=c$ and $ba = ab$.

\textbf{d.} It follows from \textbf{b} that the ring $A$ is right quasi-invariant. Therefore, $J(A)=\cap_{i\in I}M_i$, where $\{M_i\}_{i\in I}$ is the set of all ideals of the ring $A$ such that $A/M_i$ is a division ring. It is directly verified that any factor ring of a right quasi-invariant ring is right quasi-invariant. In addition, every ring $(A/M_i)[x]$ is isomorphic to a factor ring of the right quasi-invariant ring $R$. Therefore, every ring $(A/M_i)[x]$ is right quasi-invariant. It follows from \textbf{c} that every factor ring $A/M_i$ is commutative. Therefore, the factor ring $A/J(A)$ is commutative.~$\rhd$

\textbf{2.2. Lemma.} Let $A$ be a centrally essential ring.

\textbf{a.} The ring $A[x]$ is centrally essential. 

\textbf{b.} If the factor ring $A/J(A)$ is not commutative, then $A[x]$ is a centrally essential ring which is not right (left) quasi-invariant.

$\lhd$ \textbf{a.} The assertion is proved in \cite[Lemma 2.2]{MT18}. 

\textbf{b.} The assertion follows from \textbf{a} and Lemma 2.1(d).~$\rhd$

\textbf{2.3. The completion of the proof of Theorem 1.4.}

$\lhd$ \textbf{a.} Theorem 1.5 \cite{MT20c} contains an example of a centrally essential ring $A$ such that the factor ring $A/J(A)$ is not a $PI$ ring. In particular, the ring $A/J(A)$ is not commutative. Now the assertion follows from Lemma 2.2(b).

\textbf{b.} Let $\mathbb{F}$ be a field and let $\mathcal{A}$ be the subalgebra of the $\mathbb{F}$-algebra of all matrices of order $7\times 7$ consisting of matrices $A$ of the form 
$$
\left(\begin{matrix}
\alpha & a & b & c & d & e & f\\
0 & \alpha & 0 & b & 0 & 0 & d\\
0 & 0 & \alpha & 0 & 0 & 0 & e\\
0 & 0 & 0 & \alpha & 0 & 0 & 0\\
0 & 0 & 0 & 0 & \alpha & 0 & a\\
0 & 0 & 0 & 0 & 0 & \alpha & b\\
0 & 0 & 0 & 0 & 0 & 0 & \alpha\\
\end{matrix}\right).
$$
Let for $A'\in \mathcal{A}$, we have $a' = a + 1$ and the remaining components coincide with the corresponding components of the matrix $A$. Then $AA'\neq A'A$ if $a\neq 0$ and $b\neq 0$. Thus, the algebra $\mathcal{A}$ is not commutative. It is directly verified that $Z(\mathcal{A})$ consists of matrices $C$ of the form
$$
C = 
\left(\begin{matrix}
\alpha & 0 & 0 & c & d & e & f\\
0 & \alpha & 0 & 0 & 0 & 0 & d\\
0 & 0 & \alpha & 0 & 0 & 0 & e\\
0 & 0 & 0 & \alpha & 0 & 0 & 0\\
0 & 0 & 0 & 0 & \alpha & 0 & 0 \\
0 & 0 & 0 & 0 & 0 & \alpha & 0 \\
0 & 0 & 0 & 0 & 0 & 0 & \alpha\\
\end{matrix}\right).
$$
Let $A \in \mathcal{A}$ and let $a\neq 0$ or $b\neq 0$. We take a matrix $B\in Z(\mathcal{A})$ which has $d = a$, $e = b$ and zeros on the remaining positions. 
Then $0\neq AB\in Z(\mathcal{A})$. Consequently, $\mathcal{A}$ is a centrally essential algebra. 

We consider the right ideal $I$ of $\mathcal{A}$ consisting of matrices of the form 
$$
B = 
\left(\begin{matrix}
0 & 0 & b & 0 & 0 & 0 & f\\
0 & 0 & 0 & b & 0 & 0 & 0\\
0 & 0 & 0 & 0 & 0 & 0 & 0\\
0 & 0 & 0 & 0 & 0 & 0 & 0\\
0 & 0 & 0 & 0 & 0 & 0 & 0 \\
0 & 0 & 0 & 0 & 0 & 0 & b \\
0 & 0 & 0 & 0 & 0 & 0 & 0\\
\end{matrix}\right).
$$
It is directly verified that $I$ is not an ideal of $\mathcal{A}$. In addition, $I$ is a closed right ideal. Indeed, the ideal of  $\mathcal{A}$, which has only the element $c$ as a non-zero component, is a $\cap$-complement to $I$. 

At the same time, the closed left ideal $J$ of $\mathcal{A}$ consisting of matrices
$$
D = 
\left(\begin{matrix}
0 & a & 0 & 0 & 0 & 0 & f\\
0 & 0 & 0 & 0 & 0 & 0 & 0\\
0 & 0 & 0 & 0 & 0 & 0 & 0\\
0 & 0 & 0 & 0 & 0 & 0 & 0\\
0 & 0 & 0 & 0 & 0 & 0 & a \\
0 & 0 & 0 & 0 & 0 & 0 & 0 \\
0 & 0 & 0 & 0 & 0 & 0 & 0\\
\end{matrix}\right),
$$
is not an ideal. The ideal which has only the element $c$ as a non-zero component, is a $\cap$-complement of $J$, as well.~$\rhd$

\end{document}